\documentclass[12pt]{amsart}
\usepackage{latexsym,cite,amsthm,amsmath,amssymb}

\usepackage{amsmath}
\usepackage{amsthm}
\date{}

\def\a{\alpha}

\def\d{\delta}
\def\s{\sigma}

\def\0{\bar 0}
\def\1{\bar 1}

\def\Cni{\left(\begin{array}{c} 2n\\ i \end{array}\right)}
\def\Cnn{\left(\begin{array}{c} 2n\\ n \end{array}\right)}

\def\ctd{\hfill$\Box$}

\def\bes{\begin{eqnarray*}}
\def\ees{\end{eqnarray*}}
\def\bee{\begin{eqnarray}}
\def\eee{\end{eqnarray}}
\def\ˆ{^}

\def\la{{\langle}}
\def\ra{{\rangle}}
\def\Proof{{\it Proof. }}

\newtheorem{Th}{Theorem}[section]
\newtheorem{Lem}{Lemma}[section]

\newtheorem{Prop}{Proposition}[section]

\theoremstyle{remark}

\begin{document}
\title{New elements in the center of free alternative algebra}
\vspace {5 mm}

\author[Shestakov]{Ivan Shestakov}\address{Instituto de Matem\'atica e Estat\'istica,
Universidade de S\~ao Paulo,
S\~ao Paulo, Brazil
{\tiny and}
Novosibirsk State University, Novosibirsk, Russia} 
\email{shestak@ime.usp.br}

\author[Sverchkov]{Sergei Sverchkov}
\address{Novosibirsk State University\\  Novosibirsk \\  Russia}
\email{sverchkovsr@gmail.com}

\thanks{The  first author was supported by the Mathematical Center in Akademgorodok, agreement with the Ministry of Science and Higher Education of the Russian Federation No. 075-15-2019-1675.
He was also partially supported by FAPESP, Proc.\,  2018/23690-6 and CNPq, Proc.\,304313/2019-0 of Brazil.
 The second author was partially supported by FAPESP, Proc.\,2018/03717-7 of Brazil}



  \vspace {5 mm}

\maketitle
\begin{abstract}
A new series of central elements is found in the free alternative algebra.  More exactly,  let $Alt[X]$ and $SMalc[X]\subset Alt[X]$ be the free alternative algebra and the free special Malcev algebra over a field of characteristic 0 on a set of free generators $X$,  and let $f(x,y,x_1,\ldots,x_n)\in SMalc[X]$  be a multilinear element which is trivial in the free associative algebra.  Then the element $u_n=u_n(x,x_1,\ldots,x_n)=f(x\ˆ2,x,x_1,\ldots,x_n)-f(x,x\ˆ2,x_1,\ldots,x_n)$ lies in the center of the  algebra $Alt[X]$.  The elements $u_n(x,x_1,\ldots,x_n)$ are uniquely defined up to a scalar for a given $n$,  and
they are skew-symmetric on the variables $x_1,\ldots,x_n$.  Moreover, $u_n=0$ for $n=4m+2,\,4m+3$.  and  $u_n\neq 0$  for $n=4m,4m+1$.  The ideals generated by the elements  $u_{4m},\,u_{4m+1}$  lie in the associative center of the algebra $Alt[X]$ and have trivial multiplication.  
\end{abstract}

\section{An Introduction}
\hspace{\parindent}

In the study of nonassociative algebras,  various special functions play an important role, such as Kleinfeld's function in alternative algebras \cite{ZSSS}, Zelmanov's {\em tetrad-eators} in Jordan algebras \cite{Zel}, Filippov's functions in Malcev algebras \cite{F}, etc. 
In particulr, the central functions with the values in the center or  in the associative center of an algebra are very important. To find such functions, one has to investigate  the centers of free algebras.

The centers of the free alternative algebras were studied in a series of papers \cite{Dor, F, Shel, Sh-01, ShZh4, ShZh5, Sv}, where certain central elements were found. However, it is still unknown what  sets  generate the centers as  verbal subalgebras (or $T$-subalgebras).
We do not even know whether there exist  finite sets of such generators neither for associative center nor for the center.   
 In this connection, any new independent central element (that is, not lying in the $T$-subalgebra generated by the known central elements), has a special interest.

We construct a new series of central elements in the free alternative algebra.  More exactly,  let $Alt[X]$ and $SMalc[X]\subset Alt[X]$ be the free alternative algebra and the free special Malcev algebra over a field of characteristic 0 on a set of free generators $X$,  and let $f(x,y,x_1,\ldots,x_n)\in SMalc[X]$  be a multilinear element which is trivial in the free associative algebra.  Then the element $u_n=u_n(x,x_1,\ldots,x_n)=f(x\ˆ2,x,x_1,\ldots,x_n)-f(x,x\ˆ2,x_1,\ldots,x_n)$ lies in the center of the  algebra $Alt[X]$.  We prove further that the elements $u_n(x,x_1,\ldots,x_n)$ are uniquely defined up to a scalar for a given $n$,  and they are skew-symmetric on the variables $x_1,\ldots,x_n$.  Moreover, $u_n=0$ for $n=4m+2,\,4m+3$.  and  $u_n\neq 0$  for $n=4m,4m+1$.  Besides, the functions $u_n$ are zero in any alternative algebra on $n-1$ generators, which  gives one more proof of the strong inclusions 
\bes
Alt_{4m-1}\subsetneqq Alt_{4m+1};\ \ Alt_{4m}\subsetneqq Alt_{4m+2},
\ees
where $Alt_n$ means the variety of algebras generated by the free alternative algebra on $n$ generators (remind that the strong inclusions $Alt_n\subsetneqq Alt_{n+1}$ were proved in \cite{F1} for all $n$).

 The ideals generated by the elements  $u_{4m},\,u_{4m+1}$  lie in the associative center of the algebra $Alt[X]$ and have trivial multiplication.  

It remains an open question whether the verbal subalgebra $T={ alg\,_T}\{u_n,\, n=4,5,\ldots\}$ generated by the elements $u_n$ is finitely generated. We prove only that that  $T={ alg\,_T}\{u_{4m}\,| m=1,2,\cdots\}$.

\section{$\delta$-operation and $\delta$-identities.}

\hspace{\parindent}
All algebras in this section are considered over a field $F$ of characteristic 0.
We will denote by $Alt$ the variety of alternative algebras. Let $Alt[X ]$  and
$SMalc[X ]$ be the free alternative and the free special Malcev algebras with set of
generators $X=\{x_1,\ldots,,x_n\ldots\}$. Remind that the algebra $SMalc[X]$ is generated by the set $X$ as an algebra over $F$ with respect to the product $[x,y]=xy-yx$, where $xy$ is the product in $Alt[X]$.  
The elements of $SMalc[X ]$ are called {\it Malcev polynomials}.
We will denote by $x\circ y= xy+yx$ the {\em Jordan product} in the algebra $Alt[X ]$.
Let $S\subset Alt[X ]$. We will denote by $Vect_F(S)$ and $I_{Alt X}(S)$ the vector space and the
ideal of $Alt[X ]$ generated by the set $S$. The notation $I\vartriangleleft Alt[X ]$ means that $I$ is an ideal of
$Alt[X ]$. 
We will omit brackets for  left-normed products in nonassociative words. For example, $m=x_1\cdot(x_1x_2)\cdot(x_1x_2)\cdot x_1$ means
$((x_1\cdot(x_1x_2))\cdot(x_1x_2))\cdot x_1$. Set also $D(a,b,c)=(a\circ b)\circ c - (a\circ c)\circ b,\ J(a,b,c)=[[a,b],c]+[[b,c],a]+[[c,a],b]$.

Let $N = N(Alt[X ]), Z =Z(Alt[X ])$ denote the associative center and the center
of $Alt[X ]:\ N =\{a\in Alt[X ]\,|\,(a,Alt[X], Alt[X ])=0\},\ Z =\{a\in N\,|\, [a, Alt[X ]]=0 \}$.
We shall call every ideal contained in the associative center $N$ a {\em nuclear ideal}.
The ideal $I \neq 0$ of $Alt[X ]$ is called {\em trivial} if  $I^2= 0$.

This section is devoted to the construction of an infinite series of elements of $Z$.

Let $f(x,y,\ldots)\in Alt[X],\ x,y\in X,\ \deg_x (f )= \deg_y (f )=1$. For a given element $a\in Alt[X]$ denote
 $\d_af(x,y,\ldots)=f(a^2,a,\ldots)-f(a,a^2,\ldots)$. 
For briefness, we will use the symbol $\d f$ to denote $\d_af(x,y,\ldots)=\d f(a^2,\ldots)$. In particular, 
$\d f =\d g$ means that
\[
f(a^2,a,\ldots)- f(a,a^2,\ldots) = g(a^2,a,\ldots)- g(a,a^2,\ldots).
\]
For example,
\bes
\d[(a^2, x, y), (a, z, t)] &=&[(a^2, x, y), (a, z, t)] -[(a, x, y), (a^2, z, t)],\\
\d[(a^2, x, y), (a, x, y)] &=&2[(a^2, x, y), (a, x, y)].
\ees
Let $S\subset Alt[X ]$. We will denote by $\d S$ the set $\{\d f \,|\, f\in S,\, \deg_x(f)= \deg_y (f )=1,\ x, y \in X\}$.

The following well-known alternative identities may be found in \cite {ZSSS}.  We will use the same references for their  linearizations.
\bee
&([a,b],b,c) = (a,b,[c,b])=[b, (a,b, c)],\label{id_01} &\\ 
&J (a,b,c) = 6(a,b,c),&  \label{id_02}  \\
&[a^2,b]=[a,a\circ b] = [a,b]\circ a,&\label{id_03}\\ 
&(a^2,b,c) = (a,b,c) \circ a = (a,b \circ a,c),&\label{id_04}\\
&(a,b,c)\circ [a,b] = (a,b,c \circ [a,b]) = 0, &\label{id_05},\\
& ([a,b]^2 ,b,c) = ([[a,b]^2 ,b],c,d) = 0,&\label{id_06} \\
&2[J (x, y, z), t] =J ([x, y], z, t) + J ([ y, z], x, t) + J ([z, x], y, t),&  \label{id_07}\\
&J ([a,b], x, y) =[J (a, x, y),b]-[J (b, x, y), a]- 2J(a,b,[x, y]),&\label{id_08}\\
&D(a,b,c) = 2(a,b,c) +[a,[b,c]],&\label{id_09}\\
&[(a, x, y)\ˆ2,a] = 0.&\label{id_010}
\eee
An identity of type $\d f = 0$ is called a {\em $\d$-identity}. Here are some obvious examples of $\d$-identities:
\bee
\d [a^2,[a,x ]]&=&[a^2,[a ,x ]]- [a,[a^2,x ]]\nonumber \\
&=&[a ,a\circ[a,x]-[a\ˆ2,x]]=0,\label{id_011}\\
\d[a^2,(a,x,y)]&=& [a^2,(a,x,y)] -[a,(a^2,x,y)]\nonumber\\
&=&[a,a\circ(a,x,y)]- [a,a\circ(a,x,y)]=0, \label{id_012}\\
\d(a^2,(a,x,y),z)&=&(a^2,(a,x,y),z)-(a,(a^2,x,y),z)\nonumber\\
&=&(a,a\circ(a,x,y),z)-(a,a\circ(a,x,y),z)=0.\label{id_013}
\eee

In this Section we describe all $\d$-identities $\delta f$ for  $f\in SMalc[X]$. We have divided our
description into sequences of lemmas and propositions. 
\begin{Prop}\label{prop_1}
 In the variety $Alt$, the following identities are valid:
\bee
 &\d[(a^2, x, y), (a, x, y)] = 0,&\label{id_3}\\
&\d (([a^2, x],a, y)) = 0,& \label{id_4}\\
&\d([a^2, x]\circ(a, y, z)) = 0,&\label{id_5}\\
&\d((a^2, x,[[a, z], y])) = 0,&\label{id_6}\\
&\d((a^2, x, ([a, z], x, y)))= \d((a^2, x, ([a, x], z, y)))= 0,&\label{id_7}\\ 
&\d ((a^2, x,[a, z]\circ y)) = 0,&\label{id_8}\\
&\d(([a^2, x]\circ z)\circ(a, x, y)) = \d(([a^2, z]\circ x)\circ(a, x, y))=0,&\label{id_9}\\
&\d(b, (a^2, x, y), (a, x, y))=0.&\label{id_10}
\eee
\end{Prop}
\Proof
Identity \eqref{id_3}:
\bes
&[(a^2, x, y),(a, x, y)] -[(a, x, y),(a^2, x, y)]&\\
 &=[(a, x, y)\circ a, (a, x, y)]-[(a, x, y),a\circ(a, x, y)]=2[a,(a,x,y)^2]=0.&
\ees
Identity \eqref{id_4}:
\[
 ([a^2,x],a,y)= (a \circ[a, x],a, y)= ([a, x], a^2, y).
\]
Identity \eqref{id_5}: It is obvious that $[a^2, x]\circ(a, y, z),\, [a, x]\circ(a^2, y, z)$ are skew-symmetric
on $x, y, z$ . Therefore
\bes
[a^2,x]\circ (a,y,z) &=&[a,a\circ x ]\circ (a,y,z)=-[a,y]\circ (a,a\circ x,z)\\
&=&-[a,y]\circ (a^2,x,z)=[a,x]\circ (a^2,y,z).
\ees
Identity \eqref{id_6}: We have
\[
(a^2, x,[[a, z], y]) =(a, x,a \circ[[a, z], y])= (a, x,[[a^2, z], y])-(a, x,[a, y]\circ[a, z]).
\]
We have $(a, x,[a, y]\circ[a, z]) =0$ by linearization of  identity $([a,b]^2 ,b,c)=0$.
Hence $(a^2, x,[[a, z], y]) = (a, x,[[a^2, z], y])$, proving \eqref{id_6}.\\[0,5mm]
Identities \eqref{id_7}: We have
\bes
(a^2,x,([a,z],x,y))&=& (a,x,a\circ ([a,z],x,y))\\
&=&(a,x,([a^2,z],x,y)-[a,z]\circ(a,x,y))\\
&=&(a,x,([a^2,z],x,y))+(a,z,[a,x]\circ(a,x,y))\\
&\stackrel{\eqref{id_05}}=&(a,x,([a^2,z],x,y)).
\ees
Furthermore,
\bes
\d(a^2,x,([a,x],z,y))&\stackrel{\eqref{id_01}}=&\d(a^2,x,-([a,z],x,y)+(a,z,[y,x])+(a,x,[y,z]))\\
&\stackrel{\eqref{id_013}}=&-\d(a^2,x,([a,z],x,y))=0.
\ees
Identity \eqref{id_8}: We have
\bes
\d(a^2, x,[a, z]\circ y)=-(\d(a^2, x,[a, z])\circ y) +\d((a^2, x, y)\circ[a, z])\stackrel{\eqref{id_4},\eqref{id_5}}=0.
\ees
Identity \eqref{id_9}: We have
\bes
&&\d(([a^2,x]\circ z)\circ(a,x,y))=\d(a,x, ([a^2,x]\circ z)\circ y)-\d(a,x,[a^2,x]\circ z)\circ y\\
&\stackrel{\eqref{id_8}}=&\d(a,x, ([a^2,x]\circ z)\circ y)=\d(a,x,([a^2,x\circ z]-[a^2,z]\circ x)\circ y)\\
&\stackrel{\eqref{id_8}}=&-\d(a,x,([a^2,z]\circ x)\circ y)=-\d(a,x,D([a^2,z],x,y))-\d(a,x,([a^2,z]\circ y)\circ x)\\
&=&-\d(a,x,D([a^2,z],x,y))-\d(a,x^2,[a^2,z]\circ y)\stackrel{\eqref{id_8}}=-\d(a,x,D([a^2,z],x,y))\\
&\stackrel{\eqref{id_09}}=&-\d(a,x,2([a^2,z],x,y)+[[a^2,z],[x,y]])\stackrel{\eqref{id_6},\eqref{id_7}}=0.
\ees
Identities \eqref{id_10}: We have
\bes
2(b,(a^2,x,y),(a,x,y))&=&D(b,(a^2,x,y),(a,x,y))-[b,[(a^2,x,y),(a,x,y)]]\\
&\stackrel{\eqref{id_3}}=&D(b,(a^2,x,y),(a,x,y))=\d(((a^2,x,y)\circ b)\circ(a,x,y)). 
\ees
Hence
\bes
12(b,(a^2,x,y),(a,x,y))=\d((J(a^2,x,y)\circ b)\circ(a,x,y))\\
=\d((([[a^2,x],y]+[[x,y],a^2]+[[y,a^2],x])\circ b)\circ(a,x,y)). 
\ees
Therefore, it suffices to prove that 
\[
\d (([[a^2, x], y]\circ b)\circ(a, x, y)) = \d(([[x, y],a^2]\circ b)\circ(a, x, y)) = 0.
\]
For the first identity we have
\bes
\d (([[a^2, x], y]\circ b)\circ(a, x, y))& =&\d (([[a^2, x]\circ b, y]-[b,y]\circ[a^2,x])\circ(a, x, y))\\
& \stackrel{\eqref{id_9}}=&\d (([[a^2, x]\circ b, y])\circ(a, x, y))\\
& \stackrel{\eqref{id_05}}=&-\d (([x, y])\circ(a, [a^2, x]\circ b, y))\\
&=&-[x, y]\circ\d(a, [a^2, x]\circ b, y)\stackrel{\eqref{id_8}}=0.
\ees
For the second identity we have
\bes
&\d(([[x, y],a^2]\circ b)\circ(a, x, y))=\d(([[x, y]\circ b,a^2]-[b,a^2]\circ[x,y])\circ(a, x, y))&\\
&\stackrel{\eqref{id_5}}=-\d(([b,a^2]\circ[x,y])\circ(a, x, y))\stackrel{\eqref{id_9}}=\d(([b,a^2]\circ y)\circ(a, x, [x,y]))&\\
&=\d([b,a^2]\circ y)\circ[(a,x,y),x]&\\
& \stackrel{\eqref{id_03}}=\d([(a,x,y)\circ([b,a^2]\circ y),x]-(a,x,y)\circ[[b,a^2]\circ y,x])&\\
&=[\d((a,x,y)\circ([b,a^2]\circ y)),x]+\d(a,x,[b,a^2]\circ y)\circ[y,x]\stackrel{\eqref{id_8},\eqref{id_9}}=0.&
\ees
The proposition is proved.
\ctd

\smallskip

We will denote by $J=JSMalc[X ]$ and by $D=D(Alt[X ])$ the ideal of $SMalc[X]$ generated by all 
the jacobians and  the ideal of $Alt[X ]$ generated by all the associators, respectively.
It is easy to check by \eqref{id_07} that
\bee\label{id_11}
J=Vect_F\{J(a,b,x)\, |\, a,b\in SMalc[X],x\in X\}.
\eee

The largest ideal of $Alt[X ]$ contained in the center $N$  is called the {\em associative
nucleus}. We will denote it by $U =U(Alt[X ])$.  Recall the main
property of elements of $U$. Let $u\in U,\, d=d(x, y, z,\ldots)\in D$ be homogeneous of degree $\geq 1$ on x, then $d(u, y, z,\ldots)= 0$. We
have (see \cite{Sv}) for more details)
\bee\label{id_12}
\forall a,b,c\in Alt[X ], [a,[a,b]\circ[a, c]]\in U.
\eee
Let $h(a,b,\ldots,x_1,\ldots,x_n)\in Alt[X]$. Define
\bes
h_{alt(X)}&=&h(a,b,\ldots,x_1,\ldots,x_n)_{alt(X)}=\sum_{\s\in S_n}(-1)^{sgn\,\s} h(a,b,\ldots,x_{\s(1)},\ldots,x_{\s(n)}),\\
h_{sym(X)}&=&h(a,b,\ldots,x_1,\ldots,x_n)_{sym(X)}=\sum_{\s\in S_n} h(a,b,\ldots,x_{\s(1)},\ldots,x_{\s(n)}).
\ees

\begin{Lem}\label{lem_1} Let $f(x,y,z,\ldots)\in SMalc[X],\, x,y,z\in X,\, \deg_xf=\deg_yf=1$. In the algebra $Alt[X ]$, the following implications are valid:
\bee
&\deg f\leq 3\Rightarrow \d f=0,&\label{id_13}\\
&\deg f\leq 4\Rightarrow \d f\in U,&\label{id_14}\\
&\deg f\leq 5,\,f\in J \Rightarrow \d f=0,&\label{id_15}\\
&\deg f\leq 6,\,f\in J \Rightarrow \d f\in Vect_F(\d ([a^2,x_1,x_2,x_3],a,x_4)_{alt(X)}),&\label{id_16}\\
&\d([a,\underbrace{x,\ldots,x}_{\text{n}},a^2])\in U, \, n\geq 0,&\label{id_17}\\
&\d([[a,x_1,\ldots,x_k],[a^2,x_{k+1},\ldots,x_n]])\in U+\d J.\label{id_18}&
\eee
\end{Lem}
\Proof
Implication \eqref{id_13} follows from \eqref{id_011}.
For \eqref{id_14} it suffices to show that $\d([[[a, z], t], a^2])\in U$, by \eqref{id_13}. Assume first that $t=z$, then
\bes
\d([[[a, z], z],a^2]&=&[[[a, z], z],a^2]-[[[a^2, z], z],a]\\
&=&[[[a, z], z]\circ a,a]-[[[a^2, z], z],a]\\
&=&-[[a,z]^2,a]\stackrel{\eqref{id_12}}\in U.
\ees
Therefore,
$\d([[[a, z], t],a^2])$ is skew-symmetric on $z, t$ modulo $U$. Now,  we have modulo $U$
\[
2\d([[[a, z], t],a^2])\stackrel{\eqref{id_13}}\equiv\d([J (a, z, t),a^2])\stackrel{\eqref{id_012}}\equiv 0.
\]
\smallskip

Implication \eqref{id_15}: From \eqref{id_011},\eqref{id_012}, and \eqref{id_11} we have \eqref{id_15} for $\deg f <5$. Let $\deg f =5$,
then $f$ is a linear combination of associators $([a,b],[c,d],e),\, ([[a,b],c],d,e)$ for
$a,b,c,d,e\in X$, by \eqref{id_11}. It suffices to prove that $\d_a([x,z],[y,t],u)=\d_a ([[x,z],t],y,u)=0$,
by \eqref{id_011},\eqref{id_012}, and \eqref{id_13}.
We have
\bes
([a^2,z],[a,t],u)&=& (a\circ[a,z],[a,t],u)=([a,z],a\circ[a,t],u)+(a,[a,z]\circ[a,t],u)\\
&=&([a,z],[a^2,t],u)+(a,[a,z]\circ[a,t],u)\stackrel{\eqref{id_06}}=([a,z],[a^2,t],u).
\ees
and
\bes
([[a^2,z],t],a,u)&=& ([[a,z]\circ a,t],a,u)= ([[a,z],t]\circ a+ [a,z]\circ[a,t],a,u)\\
&=&([[a,z],t],a^2,u)+([a,z]\circ[a,t],a,u)\stackrel{\eqref{id_06}}=([[a,z],t],a^2,u).
\ees

Implication \eqref{id_16}: Similarly as in \eqref{id_15}, by \eqref{id_11} it suffices to prove the statement for
the associators $([[[a,b],c], d], e, f ), ([[a,b], c],[d, e],f)$ for $a,b,c,d,e, f \in X$.
By \eqref{id_13}-\eqref{id_15}, we need to consider only 2 cases: $\d([[[a^2,b],c], d], a, e)$, 
$\d([[a^2,b],c],[a,d],e)$. For the second case, we have
\bes
\d ([[a^2,b],c],[a,d],e)&\stackrel{(1)}=&-\d([[a^2,b],c],[e,d],a)\\
&+&\d([a,[[a^2,b],c]],d,e)+\d([e,[[a^2,b],c]],d,a)\\
&\stackrel{\eqref{id_15}}=&(\d([a,[[a^2,b],c]]),d,e)+\d([e,[[a^2,b],c]],d,a)\\
&\stackrel{\eqref{id_14}}=&-\d([[[a^2,b],c],e],d,a).
\ees

Consequently, it remains to prove the statement for element of type
$\d f$ for $f= ([[[a^2,z],t],u],v,a)$. Let us prove that $\d f$ is skew-symmetric on $z,t,u,v$.

We have
\bes
\d([[[a^2,z],t],u],u,a)=\d([u,([[a^2,z],t],u,a)])=[u,\d ([[a^2,z],t],u,a)]\stackrel{\eqref{id_15}}=0.
\ees
Furthermore,
\bes
&\d([[[a^2,z],t],t],u,a)=-\d([[[a^2,z],t],u],t,a)&\\
&=-\d([[[a^2,z],u],t] -[[a^2,z],[t,u]]-J([a^2,z],t,u),t,a)&\\
&\stackrel{\eqref{id_15}}=\d (J([a^2,z],t,u),t,a)\stackrel{(1)}=\d(-J([t,z],a^2,u)+[J(z,t,u),a^2]+[J(z,a^2,u),t],t,a)&\\
&\stackrel{\eqref{id_4},\eqref{id_15}}=\d ([J(z,a^2,u),t],t,a)=[t,\d(J(z,a^2,u),t,a)]\stackrel{\eqref{id_15}}=0.&
\ees
Finally, by the previous case and \eqref{id_15},
\bes
&\d([[[a^2,z],z],t],u,a)=-\d([[[a^2,z],t],z],u,a)&\\
&=\d([-[[a^2,t],z]-[a^2,[z,t]]-J(a^2,z,t),z],u,a)&\\
&=-\d([J(a^2,z,t),z],u,a)=-\d(J(a^2,z,[z,t]),u,a)=0.
\ees
Let $f=f (x, y, x_1, x_2, x_3, x_4), \deg f = 6$. By the above, it follows that
$\d f=\a \d([[[a^2,x_1],x_2],x_3],x_4,a)_{alt (X)}$, for $\a\in F$.

\smallskip

We will prove implication \eqref{id_17} by induction on $n$. For $n =1, 2$,  it
follows from \eqref{id_14}. Denote $[a,x]_k=[a,\underbrace{x,\ldots,x}_{\text{k}}]$. Assuming \eqref{id_17} to hold for $k<n$, we will prove it for $n$. Consider
\bes
[[a,x]_n,a^2]&=&[[[a,x]_{n-1},x]\circ a,a]=[[[a,x]_{n-1}\circ a,x]-[[a,x]_{n-2},x]\circ [a,x],a]\  \\
&\stackrel{\eqref{id_12}}\equiv &[[[[a,x]_{n-1}\circ a,x],a]+[[a,x]^2,[a,x]_{n-2}] \text{ (mod $U$)}.
\ees
By induction, for $k<n$  we have for any $b\in A$
\bee\label{id_19}
\d_{b}([[x,\underbrace{z,\ldots,z}_{\text{k}}],y])=[[b^2,\underbrace{z,\ldots,z}_{\text{k}}],b]-[[b,\underbrace{z,\ldots,z}_{\text{k}}],b^2]\in U.
\eee
Therefore, applying $\d_{[a,x]}$, we get by \eqref{id_19}
\bes
[[a,x]^2,[a,x]_{n-2}] =-[[a,x]_{n-2},[a,x]^2]\equiv -[[[a,x]^2,\underbrace{x,\ldots,x}_{\text{n-3}}],[a,x]]  \text{ (mod $U$)}.
\ees
If $n=3,$ the expression on the right side is zero, and if $n>3$ then it belongs to the ideal generated by the element $[[a,x]^2,x]$ which by \eqref{id_12} is contained in $U$. Therefore,
\bes
[[a,x]_n,a^2]\equiv [[[a,x]_{n-1}\circ a,x],a]  \text{ (mod $U$)}.
\ees
Assume that we have already proved that 
\bee\label{id_20}
[[a,x]_n,a^2]\equiv [[[a,x]_{n-k}\circ a,x]_k,a]  \text{ (mod $U$)}
\eee
for some $k\geq 1$. Consider
\bes
[[[a,x]_{n-k}\circ a,x]_k&=&[[[a,x]_{n-k-1}\circ a,x]-[a,x]_{n-k-1}\circ [a,x],x]_k\\
&=&[[a,x]_{n-k-1}\circ a,x]_{k+1}-[[a,x]_{n-k-1}\circ [a,x],x]_k.
\ees
Observe that 
\[
[[a,x]_{n-k-1}\circ [a,x],x]=[[[a,x]_{n-k-2},x]\circ[a,x],x]\in U \text{ by \eqref{id_12}}, 
\]
hence
\[
[[[a,x]_{n-k}\circ a,x]_k,a]\equiv [[[a,x]_{n-k-1}\circ a,x]_{k+1},a] \text{ (mod $U$)}.
\]
By induction, identity \eqref{id_20} is true for any $k<n$. In particular, for $k=n-1$  we have
\[
[[a,x]_n,a^2]\equiv [[[a,x]\circ a,x]_{n-1},a]=[[a^2,x]_n,a]  \text{ (mod $U$), proving \eqref{id_17}.}
\]
\smallskip

Similar arguments apply to \eqref{id_18}. First we prove by induction on $n$ that
\bee\label{id_21}
\d[a,x_1,\ldots,x_n,a^2]\in U+\d J.
\eee
For $n= 1, 2$ , the statement follows from \eqref{id_14}. Assume it holds for $n-1$ , we will
prove it for $n$ . By induction, we get
\bes
\d ([a,x_1,\ldots,[x_i,x_j],\ldots,x_n,a^2])\in U+\d J,
\ees
but
\bes
&\d ([a,x_1,\ldots,x_i,x_j,\ldots,x_n,a^2]-[a,x_1,\ldots,x_j,x_i,\ldots,x_n,a^2])&\\
&=\d([[a,x_1,\ldots,x_i,x_j,\ldots,x_n]-[a,x_1,\ldots,x_j,x_i,\ldots,x_n]&\\
&-[a,x_1,\ldots,[x_i,x_j],\ldots,x_n],a^2])+\d([a,x_1,\ldots,[x_i,x_j],\ldots,x_n,a^2])&\\
&=\d([J([a,x_1,\ldots,x_{i-1}],x_i,x_j),x_{j+1}\ldots,x_n]+[a,x_1,\ldots,[x_i,x_j],\ldots,x_n,a^2]) \in U+\d J.&
\ees
Therefore
\bes
&2\d ([a,x_1,\ldots,x_i,x_j,\ldots,x_n,a^2]\in &\\
& \d ([a,x_1,\ldots,x_i,x_j,\ldots,x_n,a^2]+[a,x_1,\ldots,x_j,x_i,\ldots,x_n,a^2]) +U+\d J.&
\ees
Consequently,
\bes 
(n!)\d ([a,x_1,\ldots,x_n,a^2])\in \d ([a,x_1,\ldots,x_n,a^2])_{sym(x)}+U+\d J\stackrel{\eqref{id_17}}=U+\d J.
\ees
 
Now we prove \eqref{id_18} by induction on $k,\, 0\leq k\leq n $. For $k= 0$ , the statement follows from \eqref{id_21}.
Let it holds for all $l \leq k -1$, then by induction 
\bes
&\d([[a,x_1,\ldots,x_k],[a^2,x_{k+1},\ldots,x_n]])=\d([[a,x_1,\ldots,x_{k-1}],[x_k,[a^2,x_{k+1},\ldots,x_n]]])&\\
&+\d([[[a,x_1,\ldots,x_{k-1}],[a^2,x_{k+1},\ldots,x_n]],x_k])&\\
&+\d J([a,x_1,\ldots,x_{k-1}],x_k,[a^2,x_{k+1},\ldots,x_n])\in U+\d J.&
\ees
The lemma is proved. 
\ctd

\smallskip

Fix $a\in SMalc[X]$ and define  $u_n=u_n(x_1,\ldots,x_n)=\d ([a,x_1,\ldots,x_{n-1}], a^2,x_n)$, and $U(n)=Vect_F\la u_k(v_1,\ldots,v_k)\,|\, k\leq n,\, v_1,\ldots,v_k\in SMalc[X]\ra$. Furthermore, denote by $I(n)$ the ideal of $Alt[X]$ generated by the set of commutators $\{[u, f]\,u\in U(n),\,f\in SMalc[X]\}$.

Observe that  by \eqref{id_15} and  \eqref{id_16}
\bee
&u_i=0 \hbox{ for } i\leq 3,\label{id_21'}&\\
&u_4 \hbox{ is  skew-symmetric on } x_1,\ldots,x_4.\label{id_21''}
\eee
Set $g_m(a)=[a,x_1,\ldots,x_m]$. For $a,b\in Alt[X]$ below the notation $a\equiv b$ means that $a-b\in I(n-1)$.
\begin{Prop}\label{Prop_2}
In the algebra $Alt[X]$, the following relations hold:
\bee
&u_n(x_1,\ldots,x_{n-2},x_{n-1},x_n)-u_n(x_1,\ldots,x_{n-1},x_{n-2},x_n)&\nonumber \\
 &\equiv 6\d((g_{n-3}(a),x_{n-2},x_{n-1}),a^2,x_n),\ n\geq 3,&\label{id_22'} \\
&u_n(x_1,\ldots,x_{n-3},x_{n-2},x_{n-1},x_n)-u_n(x_1,\ldots,x_{n-2},x_{n-3},x_{n-1},x_n)&\nonumber \\
&\equiv 6\d([(g_{n-4}(a),x_{n-3},x_{n-2}),x_{n-1}],a^2,x_n),\ n\geq 4.& \label{id_22''}
\eee
\bee 
 \text{\  \ \ \ For any multilinear }\! f\!=\!f(x,y,x_1,\ldots,x_n)\!\in \!J,\,  \d f\!\in \!U(n)\!+\!I(n\!-\!1). \label{id_23} 
\eee
\bee
\hbox{  For any $b,c\in SMalc[X], \  \d(g_{n-1}(a),a^2,[b,c])\in I(n)$}.\label{id_24}
\eee
\end{Prop}
\Proof
Relations \eqref{id_22'}, \eqref{id_22''} follow from  \eqref{id_02}.  We prove inclusion \eqref{id_23} by induction on $n$.
For $n\leq 4$, the statement follows from \eqref{id_16}. Assuming that it holds for $l<n$,  we will prove it for $n$. By induction
assumption and \eqref{id_11}, it suffices to prove the statement for the associators of two
types:
\bes
&(I)& \d(u_1(a,a^2,x_1,\ldots,x_k),u_2(x_{k+1},\ldots,x_{n-1}),x_n),\\
&(II)& \d(v_1(a,x_1,\ldots,x_k),v_2(a^2,x_{k+1},\ldots,x_{n-1}),x_n),
\ees
for Malcev multilinear monomials $u_i,v_j$. The statement for monomials of type (I) follows from \eqref{id_02}, \eqref{id_18}, 
and the induction assumption. For monomials of type (II), by the induction assumption  again, it
suffices to consider associators of the form
\bes
 \d([a,x_1,\ldots,x_k],[a^2,x_{k+1},\ldots,x_{n-1}],x_n),\ 0\leq k\leq n-1.
\ees
We prove it by induction on $k, 0\leq k \leq n-1$. For $k=0$, the statement is trivial. Let
it holds for all $l<k$, then
\bes
 &\d([a,x_1,\ldots,x_k],[a^2,x_{k+1},\ldots,x_{n-1}],x_n)\stackrel{\eqref{id_07}} = &\\
 &-\d([x_k,[a^2,x_{k+1},\ldots,x_{n-1}]],[a,x_1,\ldots,x_{k-1}],x_n)&\\
& -\d([[a^2,x_{k+1},\ldots,x_{n-1}],[a,x_1,\ldots,x_{k-1}]],x_k,x_n)&\\
&+2\d [([a,x_1,\ldots,x_{k-1}],x_k,[a^2,x_{k+1},\ldots,x_{n-1}]),x_n]\in U(n)+I(n-1),&
\ees
by \eqref{id_18} and the induction assumptions.

Finally,  by \eqref{id_08} we have 
\bes
\d(g_{n-1}(a),a^2,[b,c])&=&[\d(b,g_{n-1}(a),a^2),c]-[\d(c,g_{n-1}(a),a^2),b]\\
&-&2(b,c,\d([g_{n-1}(a),a^2])),
\ees
proving \eqref{id_24} and  the proposition.

\ctd
\begin{Th}\label{thm_1}
In the algebra $Alt[X]$, for all $n\geq 2$ , the following statements are valid:
\bee
u_n(x_1,\ldots,x_n)\in Z;\label{id_25}\\
u_n(x_1,\ldots,x_n) \hbox{ is skew-symmetric on $x_1,\ldots,x_n$};\label{id_26}\\
u_n(x_1,\ldots,x_{n-1},[b,c])=0.\label{id_27}
\eee
\end{Th}
\Proof
We will denote by $A_n, B_n$  statements \eqref{id_25} and \eqref{id_26}, respectively. Set
$D_n=A_n\& B_{n+1}$. We prove the truth of  statement $D_n$ by induction on $n$. For $n\leq 3$,
the statement follows from \eqref{id_21'}, \eqref{id_21''}. Assuming that statement $D_n$ holds for $n-1$, we will
prove it for $n\geq 4$.

In order  to prove  relation \eqref{id_25}, we will need some additional identities.
First, we prove the identity
\bee\label{id_29}
u_n(x_1,\ldots,x_{n-2},b,[x_{n-1},b])-u_n(x_1,\ldots,x_{n-1},b,[x_{n-2},b])=0.
\eee
Notice that by induction assumption $D_{n-1}$,  the ideal $I(n-1)=0,$ and  the element $u_n(x_1,\ldots,x_n)$ is skew-symmetric on $x_1,\ldots,x_n$.

 Therefore,  by \eqref{id_22'} we get
\bes
&2u_n(x_1,\ldots,x_{n-2},b,[c,b])=2u_n(x_1,\ldots,[c,b],x_{n-2},b)&\\
&=6\d ((g_{n-3}(a),[c,b],x_{n-2}),a^2,b),&
\ees
and 
\bes
&u_n(x_1,\ldots,x_{n-2},b,[x_{n-1},b])-u_n(x_1,\ldots,x_{n-1},b,[x_{n-2},b])&\\
&=3\d (((g_{n-3}(a),[x_{n-1},b],x_{n-2}),a^2,b)-((g_{n-3}(a),[x_{n-2},b],x_{n-1}),a^2,b))&\\
&=6\d([(x_{n-1},b,x_{n-2}),g_{n-3}(a)],a\ˆ2,b)\in I(n-1)=0,
\ees
proving \eqref{id_29}.

Next we prove that 
\bee\label{id_30}
[u_n(x_1,\ldots,x_{n-1},b),b]=0.
\eee
By \eqref{id_22''} we have
\bes
&2[u_n(x_1,\ldots,x_{n-1},b),b]=-6[\d([(g_{n-4}(a),x_{n-3},x_{n-2}),b],a^2,x_{n-1}),b]&\\
&\stackrel{\eqref{id_01}}=6[\d([(g_{n-4}(a),x_{n-3},x_{n-2}),b],b,x_{n-1}),a^2]&\\
&-6\d([b,[(g_{n-4}(a),x_{n-3},x_{n-2}),b]],a^2,x_{n-1})&\\
&-6\d([a\ˆ2,[(g_{n-4}(a),x_{n-3},x_{n-2}),b]],b,x_{n-1}),&
\ees
and we will prove that all the three summonds are zero.

First we have
\bes
&2[\d([(g_{n-4}(a),x_{n-3},x_{n-2}),b],b,x_{n-1}),a^2]&\\
&=2[\d((g_{n-4}(a),x_{n-3},x_{n-2}),b,[x_{n-1},b]),a^2]&\\ \ 
&\stackrel{\eqref{id_07}}=\d(([(g_{n-4}(a),x_{n-3},x_{n-2}),b],[x_{n-1},b],a\ˆ2)&\\
&+([b,[x_{n-1},b]],(g_{n-4}(a),x_{n-3},x_{n-2}),a\ˆ2)&\\
&+([[x_{n-1},b],(g_{n-4}(a),x_{n-3},x_{n-2})],b,a\ˆ2))&\\
&\stackrel{\eqref{id_22'},\eqref{id_22''}} 
=-\tfrac{1}{3} u_{n}(x_1,\ldots,x_{n-2},b,[x_{n-1},b])&\\
&+\tfrac13u_{n-1}(x_1,\ldots,[b,[x_{n-1},b]])
+\tfrac{1}{3} u_{n}(x_1,\ldots,x_{n-2},b,[x_{n-1},b])&\\
&=\tfrac13u_{n-1}(x_1,\ldots,x_{n-2},[b,[x_{n-1},b]])\stackrel{\eqref{id_24}}\in I(n-1)=0.&
\ees
Furthermore, by \eqref{id_23} and \eqref{id_02}  we have,
\bes
&\d([a\ˆ2,[(g_{n-4}(a),x_{n-3},x_{n-2}),b]],b,x_{n-1})&\\
&=(\d([a\ˆ2,[(g_{n-4}(a),x_{n-3},x_{n-2}),b]]),b,x_{n-1})\in I(n-1)=0.&
\ees
Finally,  by \eqref{id_07} and by induction, we have
\bes
&2\d([b,[(g_{n-4}(a),x_{n-3},x_{n-2}),b]],a^2,x_{n-1})=\d([b,([g_{n-4}(a),x_{n-3}],x_{n-2},b)&\\
&+([x_{n-3},x_{n-2}],g_{n-4}(a),b)+([x_{n-2},g_{n-4}(a)],x_{n-3},b)],a\ˆ2,x_{n-1})&\\
&=\d([b,([g_{n-4}(a),x_{n-3}],x_{n-2},b)+([x_{n-2},g_{n-4}(a)],x_{n-3},b)],a\ˆ2,x_{n-1})&\\
&\stackrel{\eqref{id_01}}=\d(([g_{n-4}(a),x_{n-3}],[x_{n-2},,b],b)+([x_{n-2},g_{n-4}(a)],[x_{n-3},b],b),a\ˆ2,x_{n-1})&\\
&=\tfrac13 (u_n(x_1,\ldots,x_{n-3},[x_{n-2},b],b,x_{n-1})-u_n(x_1,\ldots,x_{n-2},[x_{n-3},b],b,x_{n-1}))\stackrel{\eqref{id_29}}=0,&
\ees
 proving \eqref{id_30}.

Therefore,  to prove \eqref{id_25}, it suffices to prove that
\bee\label{id_31}
[u_n(x_1,\ldots,b),c]-[u_n(x_1,\ldots,c),b]=0.
\eee
We have by \eqref{id_08}
\bes
[u_n(x_1,\ldots,x_{n-1},b),c]-[u_n(x_1,\ldots,x_{n-1},c),b]\\
=[\d(g_{n-1}(a),a\ˆ2,b),c]-[\d(g_{n-1}(a),a\ˆ2,c),b]\\
=\d([g_{n-1}(a),a\ˆ2],b,c)+2\d(g_{n-1}(a),a\ˆ2,[b,c])
\ees
By \eqref{id_21},  $\d([g_{n-1}(a),a\ˆ2])=u+\d f$, where $u\in U,\, f\in J,  \deg f=n+1$, and we have by \eqref{id_23}
\bes
\d([g_{n-1}(a),a\ˆ2],b,c)=(\d f,b,c)\in I(n-1)=0.
\ees
Furthermore, by \eqref{id_22'} and the skew-symmetry  of $u_n$,
\bes
2\d(g_{n-1}(a),a\ˆ2,[b,c])=6\d((g_{n-3}(a),x_{n-2},x_{n-1}),a\ˆ2,[b,c])\\
=-6\d((g_{n-3}(a),x_{n-2},[b,c]),a\ˆ2,x_{n-1}).
\ees
Observe that by \eqref{id_30} the left part of \eqref{id_31} is skew-symmetric on $x_1,\ldots,x_{n-1},b,c$. 
Therefore,
\bes
&6\d((g_{n-3}(a),x_{n-2},[b,c]),a\ˆ2,x_{n-1})=&\\
&=2\d((g_{n-3},x_{n-2},[b,c])+(g_{n-3},b,[c,x_{n-2}])+(g_{n-3},c,[x_{n-2},b]),a\ˆ2,x_{n-1})&\\
&\stackrel{\eqref{id_07}}=-4\d([(b,c,x_{n-2}),g_{n-3}(a)],a\ˆ2,x_{n-1})\in I(n-1)=0.&
\ees
This proves \eqref{id_25} and the statement $A_n$.  

\smallskip
To finish the proof of $D_n$,  it remains to prove that $u_{n+1}(x_1,\ldots,x_{n+1})$ is skew-symmetric on $x_1,\ldots,x_{n+1}$. We have
\bes
&\d([a,x_1,\ldots,x_n],a\ˆ2,x_{n+1})\stackrel{\eqref{id_01}}=-\d([a,x_1,\ldots,x_{n+1}],a\ˆ2,x_{n})&\\
&+\d([x_n,([a,x_1,\ldots,x_{n-1}],a\ˆ2,x_{n+1})])+\d([x_{n+1},([a,x_1,\ldots,x_{n-1}],a\ˆ2,x_{n})])&\\
&\stackrel{\eqref{id_25}}=-\d([a,x_1,\ldots,x_{n+1}],a\ˆ2,x_{n}),&
\ees
that is, $u_{n+1}$ is skew-symmetric on $x_n,x_{n+1}$.  Furthermore,
\bes
&\d([a,x_1,\ldots,x_n],a\ˆ2,x_{n+1})\stackrel{\eqref{id_01}}=-\d([a,x_1,\ldots,x_{n-1},a\ˆ2],x_n,x_{n+1})&\\
&+\d([a,x_1,\ldots,x_{n-1}],a\ˆ2,[x_{n+1},x_n])+\d([a,x_1,\ldots,x_{n-1}],x_n,[x_{n+1},a\ˆ2])&\\
&\stackrel{\eqref{id_21},\eqref{id_24},\eqref{id_25}}=\d([a,x_1,\ldots,x_{n-1}],x_n,[x_{n+1},a\ˆ2]).&
\ees
Manipulating with $x_{n-1},x_n$ as before with $x_n,x_{n+1}$, we get
\bes
&\d([a,x_1,\ldots,x_{n-1}],x_n,[x_{n+1},a\ˆ2])=-\d([a,x_1,\ldots,x_{n-2},x_{n}],x_{n-1},[x_{n+1},a\ˆ2])&\\
&+[x_n,\d([a,x_1,\ldots,x_{n-2}],x_{n-1},[x_{n+1},a\ˆ2])]+[x_{n-1},\d([a,x_1,\ldots,x_{n-2}],x_{n},[x_{n+1},a\ˆ2])]&\\
&\stackrel{\eqref{id_25}}=-\d([a,x_1,\ldots,x_{n-2},x_{n}],x_{n-1},[x_{n+1},a\ˆ2])&
\ees
Hence $u_{n+1}$ is skew-symmetric on $x_{n-1}$ and $x_n$. Repeating in the same way, we get
\bee
&\d([a,x_1,\ldots,x_n],a\ˆ2,x_{n+1})=\d([a,x_1,\ldots,x_{n-1}],[a\ˆ2,x_{n+1}],x_n)&\nonumber \\
&=\d([a,x_1,\ldots,x_{n-2}],[[a\ˆ2,x_{n+1}],x_n],x_{n-1})=\ldots&\nonumber\\
&\ldots=\d([a,x_1],[a\ˆ2,x_{n+1},\ldots,x_3],x_2)=\d(a,[a\ˆ2,x_{n+1},\ldots,x_2],x_1),&\label{id_32}
\eee
where on the $i$-th pass we have, as above,
\bes
&\d([a,x_1,\ldots,x_{i-1},x_i],[a\ˆ2,x_{n+1},\ldots,x_{i+2}],x_{i+1})&\\
&=-\d([a,x_1,\ldots,x_{i-1},x_{i+1}],[a\ˆ2,x_{n+1},\ldots,x_{i+2}],x_{i}).&
\ees
This proves the skew-symmetry of $u_{n+1}$ on its arguments. Thus $D_{n+1}$ is proved and therefore \eqref{id_25},\eqref{id_26} are proved.

Finally, identity \eqref{id_27} follows from \eqref{id_24} and \eqref{id_25}.

\ctd

\section{Further properties of functions $u_n$.}

\begin{Prop}\label{prop_3}
For any $a_1,\ldots,a_n,a\in Alt[X]$, the element $u_n(a;a_1,\ldots,a_n)=\d_a(g_{n-1}(a),a\ˆ2,a_n)$ lies in the ideal generated by the elements $u_n(a;x_{i_1},\ldots,x_{i_n})$,  where $x_j\in X$. In particular,  any $n$-generated alternative algebra satisfies the identity $u_{n+1}=0$.
\end{Prop}
\Proof 
By \eqref{id_04} we have 
\bee\label{id_33}
&&u_n(x_1,\ldots,x_{n-1},x\circ y)=\d(g_{n-1}(a),a\ˆ2,x\circ y)\nonumber \\
&=&\d(g_{n-1}(a),a\ˆ2,x)\circ y+\d(g_{n-1}(a),a\ˆ2,y)\circ x\nonumber \\
&=&u_n(x_1,\ldots,x_{n-1},x)\circ y+u_n(x_1,\ldots,x_{n-1},y)\circ x.
\eee
Now the Proposition follows from \eqref{id_26}, \eqref{id_27}.

\ctd
\begin{Prop}\label{prop_4}
$u_n(a;x_1,\ldots,x_n)=0$ for any $a\in Alt[X]$ and $n=4k+2,\,4k+3,\ k\geq 0$.
\end{Prop}
\Proof
In view of \eqref{id_32} we have
\bes
&\d([a,x_1,\ldots,x_{n-1}],a\ˆ2,x_{n})=\d(a,[a\ˆ2,x_{n},\ldots,x_2],x_1)&\\
&=-\d([a\ˆ2,x_{n},\ldots,x_2],a,x_1)=\d([a,x_{n},\ldots,x_2],a\ˆ2,x_1)&\\
&=(-1)\ˆm \d([a,x_1,\ldots,x_{n-1}],a\ˆ2,x_{n}),&
\ees
where $m=1+2+\ldots+(n-1)=\tfrac{n(n-1)}{2}$. Therefore,
\bes
u_n(x_1,\ldots,x_n)=\d([a,x_1,\ldots,x_{n-1}],a\ˆ2,x_{n})=-\d([a,x_1,\ldots,x_{n-1}],a\ˆ2,x_{n})=0
\ees
for $n=4k+2,\,4k+3$.

\ctd

In order to prove that $u_n(x;x_1,\ldots,x_n)\neq 0$ for $n=4k,\, 4k+1,\ k>0,$ 
it  suffices to show that there exist an alternative algebra $A$ and elements $e,e_1,\ldots,e_n\in A$ such that
\bes
\d([e,e_1,\ldots,e_{n-1}],e\ˆ2,e_n)\neq 0 \hbox{ for all } n=4k,\, 4k+1,\ k>0.
\ees
Observe that by \eqref{id_02} and Proposition \ref{Prop_2} the previous inequalities up to a nonzero scalar may be rewritten as
\bes
\d([(e\ˆ2,e_1,e_2),\ldots,e_{n-3},e_{n-2}),e_{n-1}],e_n,e)\neq 0, \ n=4k, \ k>0.\\
\d((e\ˆ2,e_1,e_2),\ldots,e_{n-2},e_{n-1}),e_n,e)\neq 0, \ n=4k+1,\  k>0.
\ees

Due to skew-symmetry of the left parts of these inequalities on the variables $e_1,\ldots e_n$,  it is more convinient to prove the following superversions of these inequalities in the free alternative superalgebra $SAlt[e;x]$ generated by an even generator $e$ and  an odd generator $x$ (see, for instance, \cite{Sh-03,ShZh4}):
\bee
\d(S_n(e\ˆ2,e,x))\neq 0, &\ n=4k, \ k>0.\label{id33}\\
\d(T_n(e\ˆ2,e,x))\neq 0,&\ n=4k+1,\  k>0.\label{id34}
\eee
where 
\bes
S_n(a,b,x)&=&(([(a,\underbrace{x,x),\ldots,x,x)}_{n-2},x],x,b),\\
T_n(a,b,x)&=&(((a,\underbrace{x,x),\ldots,x,x)}_{n-1},x,b).\
\ees
\begin{Prop}\label{prop_5}
The inequalities \eqref{id33}, \eqref{id34} hold in the free alternative superalgebra $SAlt[e;x]$. 
\end{Prop}
\Proof
We will use the example of an alternative  superalgebra from \cite{Sh-04} which modifies and corrects the example of an alternative algebra by Yu.\,Med\-ve\-dev \cite{Med}.  For covinience of the reader, we give below the multiplication table of this example.

\smallskip

Consider a vector space $A_n, \, n=4k+2,$ with a base
\[
x, v_i,\,v'_j,\,u_i,\,u'_j,\, U,\, V.
\]
where $i=0,\ldots,n;\  j=1,\ldots,n$.
Denote also
$
e=v_0,\  e^2=u_0,
$
and let $w\in\{u,u',v,v'\}$.

\smallskip
Define a multiplication on $A_n$ by the following rules:\\
For $i>0$
\bes
v'_i\cdot e&=&-(-1)\ˆ{i}u'_i,\\
v_i\cdot e&=&(-1)\ˆ{i}(u_i+u'_i)\\
e\cdot v'_i&=&\left\{\begin{array}{cc}
-u_i,  &\hbox{ if $i$ even,}\\
-u_i-u'_i,&\hbox{ if $i$ odd,}
\end{array}
\right.\\
e\cdot v_i&=&\left\{\begin{array}{cc}
-u'_i,  &\hbox{ if $i$ even,}\\
u_i,&\hbox{ if $i$ odd,}
\end{array}
\right.
\ees
and for $i<n$
\bes
w_i\cdot x&=&w_{i+1},\\
x\cdot v_i&=&\left\{\begin{array}{cc}
v'_{i+1} &\hbox{ if $i$ even,}\\
v_{i+1}+v’_{i+1},&\hbox{ if $i$ odd, }
\end{array}
\right.\\
x\cdot v'_i&=&\left\{\begin{array}{cc}
-v_{i+1} &\hbox{ if $i$ odd,}\\
-v_{i+1}-v'_{i+1},&\hbox{ if $i$ even,}
\end{array}
\right.\\
x\cdot u_i&=&\left\{\begin{array}{cc}
u'_{i+1} &\hbox{ if $i$ even,}\\
u_{i+1}+u’_{i+1},&\hbox{ if $i$ odd, }
\end{array}
\right.\\
x\cdot u'_i&=&\left\{\begin{array}{cc}
-u_{i+1} &\hbox{ if $i$ odd,}\\
-u_{i+1}-u'_{i+1},&\hbox{ if $i$ even,}
\end{array}
\right.\\
\ees
Set also for $k=0,\ldots,\tfrac{n}{2},$
\bes
u_{n-2k}v_{2k}&=&v_{2k+1}u_{n-2k-1}=v_{2k}u'_{n-2k}=-u'_{n-2k-1}v_{2k+1}\\
&=&u'_{n-2k}v'_{2k}=v'_{2k+1}u'_{n-2k-1}=(-1)\ˆkU,\\[1mm]
v_{2k}u_{n-2k}&=&-u_{n-2k-1}v_{2k+1}=v'_{2k+1}u_{n-2k-1}=u_{n-2k}v'_{2k}\\
&=&v'_{2k}u'_{n-2k}=-u'_{n-2k-1}v'_{2k+1}=(-1)\ˆkV,\\[1mm]
u'_{n-2k}v_{2k}&=&-u_{n-2k-1}v'_{2k+1}=v_{2k+1}u'_{n-2k-1}=v'_{2k}u_{n-2k}=(-1)\ˆ{k-1}(U+V).
\ees
All  other products are zero.

It is proved in \cite{Sh-04} that the superalgebra $A_n$ is alternative.

To prove that $\d(S_n(e\ˆ2,e,x))\neq 0$, consider its linearization
\bes
\tilde{S}_n(e,z,t,x)&=&S_n(zt-tz,x,e)+S_n(e\circ z,x,t)- S_n(e\circ t,x,z)\\
&-&S_n(e,x,zt-tz)+S_n(t,x,e\circ z)-S_n(z,x,e\circ t)
\ees
 with respect to $e$,  with odd variables $z,t$. It suffices to prove that 
 $\tilde{S}_n(e,v_1,v'_1,x)\neq 0$  in $A_n$.

We have
\bes
 \tilde{S}_n(e,v_1,v'_1,x)&=&([(\cdots(e\circ v_1,\underbrace{x,x),\ldots,x,x)}_{\text{n-4}},x],x,v'_1)\\
&-&([(\cdots(e\circ v'_1,\underbrace{x,x),\ldots,x,x)}_{\text{n-4}},x],x,v_1)\\
&+&([(\cdots(v'_1,\underbrace{x,x),\ldots,x,x)}_{\text{n-4}},x],x,e\circ v_1)\\
&-&([(\cdots(v_1,\underbrace{x,x),\ldots,x,x)}_{\text{n-4}},x],x,e\circ v'_1)\\
&=&-([u'_{n-3},x],x,v'_1)+([u_{n-3},x],x,v_1)\\
&-&([v'_{n-3},x],x,u'_1)+([v_{n-3},x],x,u_1)\\
&=&(u_{n-2}-u'_{n-2},x,v'_1)+(2u_{n-2}+u'_{n-2},x,v_1)\\
&+&(v_{n-2}-v'_{n-2},x,u'_1)+(2v_{n-2}+v'_{n-2},x,u_1)\\
&=&u_{n-1}v'_1-u'_{n-1}v'_1+u_{n-2}v_2+u_{n-2}v'_2-u'_{n-2}v_2-u'_{n-2}v'_2\\
&+&2u_{n-1}v_1+u'_{n-1}v_1-2u_{n-2}v_2-2u_{n-2}v'_2-u'_{n-2}v_2-u'_{n-2}v'_2\\
&+&v_{n-1}u'_1-v'_{n-1}u'_1+v_{n-2}u_2+v_{n-2}u'_2-v'_{n-2}u_2-v'_{n-2}u'_2\\
&+&2v_{n-1}u_1+v'_{n-1}u_1-2v_{n-2}u_2-2v_{n-2}u'_2-v'_{n-2}u_2-v'_{n-2}u'_2
\ees
\bes
&=&U+V+V-U-V-U-V+U\\
&-&2V-U+2U+2V-U-V+U\\
&-&U-V-U+V+U+U+V-V\\
&+&2U+V-2V-2U+U+V-V=2U-2V\neq 0.
\ees

To prove that $\d(T_n(e\ˆ2,e,x))\neq 0$, consider its partial linearization
\bes
{T'}_n(e,a,z,x)&=&T_n(e\circ a,z,x)-T_n(e\circ z,a,x)- T_n(a\circ z,e,x)\\
&+&T_n(z,e\circ a,x)-T_n(a,e\circ z,x)-T_n(e,a\circ z,x)
\ees
 with respect to $e$,  with even variable $a$ and  odd variable $z$.  Let, furthermore, 
 $T''_n(e,a,z,t,x)=t\tfrac{\partial}{\partial x}({T'}_n(e,a,z,x)$ with odd variable $t$.
It suffices to prove that $T''_n(e,v_4,x,v_1,x)\neq 0$ in $A_n$.

It is easy to see that $(A_nA_n)A_n=A_n(A_nA_n)=0$, therefore
\bes
T''_n(e,v_4,x,v_1,x)&=&(\cdots(v_4\circ e,\underbrace{x,x),\cdots,x,x)}_{4n},v_1,x)\\
&-&(\cdots(v_4\circ x,\underbrace{x,x),\cdots,x,x)}_{4n},v_1,e)\\
&-&(\cdots(x\circ e,\underbrace{x,x),\cdots,x,x)}_{4n},v_1,v_4)\\
&+&(\cdots(x,\underbrace{x,x),\cdots,x,x)}_{4n},v_1,v_4\circ e)\\
&-&(\cdots(e,\underbrace{x,x),\cdots,x,x)}_{4n},v_1,v_4\circ x)\\
&-&(\cdots(v_4,\underbrace{x,x),\cdots,x,x)}_{4n},v_1,x\circ e)\\
&+&2(\cdots(x,v_1,x)\underbrace{x,x),\cdots,x,x)}_{4n-2},x,v_4\circ e)\\
&-&2(\cdots(e,v_1,x)\underbrace{x,x),\cdots,x,x)}_{4n-2},x,v_4\circ x)
\ees
\bes
&=&(\cdots(v_4\circ e,\underbrace{x,x),\cdots,x,x)}_{4n},v_1,x)\\
&-&(\cdots(v_4\circ x,\underbrace{x,x),\cdots,x,x)}_{4n},v_1,e)\\
&+&2(\cdots(x,v_1,x)\underbrace{x,x),\cdots,x,x)}_{4n-2},x,v_4\circ e)\\
&-&2(\cdots(e,v_1,x)\underbrace{x,x),\cdots,x,x)}_{4n-2},x,v_4\circ x)\\
&=&-u_{n-2}v_2-(v_{n-1}+v'_{n-1})(u_1+u'_1)+2v_{n-4}u_4-2v_{n-5}u'_5\\
&-&2(u_{n-5}+u'_{n-5})(v_5+v'_5)+2(u_{n-6}+u'_{n-6})v'_6\\
&=&U-U-V+U+V-U-2V-2U+2V\\
&-&2U-2V+2U+2V-2V-2U=-4U-2V\neq 0.
\ees

\ctd

Let $Alt_n$ denote the variety generated by the free alternative algebra with $n$ generators. It is known that  $Alt_n\subsetneqq Alt_{n+1}$ for all $n\geq 1$ (see \cite{F1}).
Propositions \ref{prop_3} and \ref{prop_5} give another proof that 
\bes
Alt_{4m-1}\subsetneqq Alt_{4m+1},\ \ Alt_{4m}\subsetneqq Alt_{4m+2}
\ees
for all $m\geq 1$.

\begin{Th}\label{thm_2}
The element $u_n(x;x_1,\ldots,x_n)$ generates a nonzero  trivial nuclear ideal in $Alt[X]$ for any $n=4m, 4m+1,\ m\geq 1$. 
\end{Th}
\Proof
Let $I=I_{Alt[X]}(u_n)$, in view of Proposition \ref{prop_5} it suffices to prove that $I\subset N$ and $I^2=0$.
Prove first that 
\bee\label{id_34}
u_n(x_1,\ldots,x_{n-1},b)\circ[a,b]=0.
\eee
In fact, by  \eqref{id_33} we have
\bes
u_n(x_1,\ldots,x_{n-1},b)\circ[a,b]&=&u_n(x_1,\ldots,x_{n-1},b\circ [a,b]) -u_n(x_1,\ldots,x_{n-1},[a,b])\circ b\\
&\stackrel{\eqref{id_27}}=&u_n(x_1,\ldots,x_{n-1},[a,b^2])=0. 
\ees
Now by \eqref{id_34},
\bes
u_n(x_1,\ldots,x_{n})\circ [[a,b],c]=-u_n(x_1,\ldots,x_{n-1},[a,b])  \circ [x_n,c]=0,
\ees
and therefore $u_n(x_1,\ldots,x_{n})\circ (a,b,c)=0$, proving that $u_n(x_1,\ldots,x_{n})\in U$. 
Finally,
\bes
2u_n^2&=&u_n(x_1,\ldots,x_{n})\circ u_n\stackrel{\eqref{id_33}}=u_n(x_1,\ldots,x_{n}\circ u_n)-u_n(x_1,\ldots,u_{n})\circ x_n\\
&=&u_n(x_1,\ldots,u_{n}(x_1,\ldots,x_n^2))-u_n(x_1,\ldots,u_{n})\circ x_n\stackrel{\eqref{id_27}}=0
\ees
since $u_n\subset [Alt[X],Alt[X]]$. 

Since by Theorem \ref{thm_1} $u_n\in Z(Alt[X])$, this finishes the proof.

\ctd

\section{ Some known results and open questions}

\hspace{\parindent} We resume here the known results and open questions on the structure of the center $Z$ of the free alternative algebra $Alt[X]$.

\smallskip

The first example of nonzero elements in $Z$ was found independently by Dorofeev \cite{Dor} and Shelipov \cite{Shel}
\bes
\hbox{ \bf Dorofeev and Shelipov (1973):  }
 [(x,y,z),t]^4\in Z.
\ees
The next example was found by Shestakov \cite{Sh-01}
\bes
\hbox{ \bf Shestakov (1976):  }
(x,y,z)^4\in Z.
\ees
Filippov in \cite{F} found an element from $Z$ of  degree 7. Till now, it remains an element of the smallest known degree from $Z$. 
In this connection, we mention\\[1mm]
{\bf  The Filippov's Conjecture:} {\it There are no nonzero elements of degree less then 7 in $Z$}.\\[1mm]
We prefer to write the Filippov's central element and the other new central elements in the superized form.
Let $SAlt[A;X]$ denotes the free alternative superalgebra on a set of  even elements $A$ and a set of odd elements $X$.  
\bes
\hbox{ \bf FIlippov (1999)}&:&   \hbox{ \it For any  $a\in A$ and $x\in X$,}\\
Fil(2a,5x)&=&([x,[xx,x]]_s\circ a-[x,[xx,x]\circ a]_s,a,x)\\
&-&((x\circ [xx,a])\circ a-x\circ([xx,a]\circ a),x,x)\in Z(SAlt[A;X]).
\ees
One can obtain the corresponding  element  from $Z$ by skew-symmetrization (see \cite{Sh-03, ShZh4}):
 \bes
Fil(a,x_1,\ldots,x_5)&=&(([x_1,[x_2x_3,x_4]]\circ a-[x_1,[x_2x_3,x_4]\circ a],a,x_5)\\
&-&((x_1\circ [x_2x_3,a])\circ a-x_1([x_2x_3,a]\circ a),x_4,x_5))_{alt(X)}\in Z.
\ees
Another central element of degree 7 was found by Hentzel and Peresi \cite{HenPer} via computer calculations. In \cite{ShZh5}, it was proved that the Hentzel-Peresi element is just the first member of an infinite series of skew-symmetric central elements.
Define, for an odd element $x\in SAlt[A;X]$, $x^{[1]}=x,\, x^{[i+1]}=[x^{[i]},x]_s$.\\

{\bf Shestakov-Zhukavets (2006):}  The element  $z_n=[x^{[n]},xx]\in Z(SAlt[A;X]) $ for any  $n\geq 5.$\\

The Hentzel-Peresi element  is just a skew-symmetrization of $z_5$. It was proved in \cite{ShZh5} that all $z_n$ for $n\geq 5$ lie in the verbal subsuperalgebra generated by $z_5$.

Finally, our first nonzero central element $u_4(a,x_1,\ldots,x_4)$ provides one more central element of degree 7.  Its super-version can be written as 
\bes
u_4(3a,4x)=([[aa,x]\circ_s x,x]\circ_s x,x],a,x)-([[a,x]\circ_s x,x]\circ_s x,x],aa,x).
\ees
Therefore, we have the three elements from $Z(Alt[A;X])$ of degree 7:\\
\bee
\label{3elements}
Fil(2a,5x),\, z_5, \, u_4(3a,4x).
\eee
We check by computer that these three elements are independent, that is, no one of them lies in the verbal subalgebra generated by the other two elements.  In this connection, we ask the question:\\
{\bf Question:} {\it Is it true that the skew-symmetrizations of elements \eqref{3elements} generate the verbal subspace of the degree 7 elements from $Z$?}

It remains an open question whether all central elements $u_n$ are independent? We can prove only the following
\begin{Prop}\label{prop_6}
The  central element $u_{4m+1}$ lie in the verbal subsuperalgebra generated by $u_{4m}$.
\end{Prop}
\Proof
Denote, for a homogeneous $b$ and odd $x$, $b_0=b,\ b_{i+1}=[b_i,x]_s$, then $u_n=((a^2)_{n-1},a,x)-(a_{n-1},a^2,x)$. The linearization of $u_{4m}$ on $a$  with an odd element $y$ lies in the verbal subsuperalgebra generated by $u_{4m}$ and has a form
\bes
u_{4m}(a,y,x)&=&((a^2)_{4m-1},y,x)-((a\circ y)_{4m-1},a,x)\\
&-&(a_{4m-1},a\circ y,x)+(y_{4m-1},a^2,x).
\ees
Substituting $y=a_1$, we have
\bes
u_{4m}(a,a_1,x)&=&((a^2)_{4m-1},a_1,x)-(((a^2)_1)_{4m-1},a,x)\\
&-&(a_{4m-1},(a^2)_1,x)+((a_1)_{4m-1},a^2,x)\\
&=&((a^2)_{4m-1},a_1,x)-(((a^2)_{4m},a,x)\\
&-&(a_{4m-1},(a^2)_1,x)+(a_{4m},a^2,x)\\
&=&-u_{4m+1}+((a^2)_{4m-1},a_1,x)-(a_{4m-1},(a^2)_1,x).
\ees
It follows from  \eqref{id_32} that
\bee
u_{n+1}=(-1)^{\tfrac{i(i+1)}{2}}(((a^2)_{n-i},a_i,x)-(a_{n-i},(a^2)_i,x)),
\eee
hence $((a^2)_{4m-1},a_1,x)-(a_{4m-1},(a^2)_1,x)=-u_{4m+1}(a,x)$, and finally $2u_{4m+1}=-u_{4m}(a,a_1,x)$.

\ctd

{\bf Question:} {\it Is it true that all central elements $u_{n}$ lie in the verbal subsuperalgebra generated by $u_{4}$?}

\end{document}